# Optimized electrified meeting-point-based feeder bus services with capacitated charging stations and partial recharges

Tai-Yu Ma[a*], Yumeng Fang[a], Richard D. Connors[b], Francesco Viti[b], Haruko Nakao[b]

[a]*Luxembourg Institute of Socio-Economic Research, 11 Porte des Sciences, Esch-sur-Alzette, L-4366, Luxembourg*
[b]*Department of Engineering, University of Luxembourg, Esch-sur-Alzette, L-4366, Luxembourg*

**Abstract**

Meeting-point-based feeder services using EVs have good potential to achieve an efficient and clean on-demand mobility service. However, customer-to-meeting-point, vehicle routing, and charging scheduling need to be jointly optimized to achieve the best system performance. To this aim, we assess the effect of different system parameters and configure them based on our previously developed hybrid metaheuristic algorithm. A set of test instances based on morning peak hour commuting scenarios between the cities of Arlon and Luxembourg are used to evaluate the impact of the set parameters on the optimal solutions. The experimental results suggest that higher meeting point availability can achieve better system performance. By jointly configuring different system parameters, the overall system performance can be significantly improved (-10.8% total kilometers traveled by vehicles compared to the benchmark) to serve all requests. Our experimental results show that the meeting-point-based system can reduce up to 70.2% the fleet size, 6.4% the in-vehicle travel time and 49.4% the kilometers traveled when compared to a traditional door-to-door dial-a-ride system.

## 1. Introduction

The climate crisis urges electrification in the transportation sector, including demand-responsive transportation (DRT) systems which introduce additional constraints regarding electric charging operations. Besides, the "meeting-point **(MP)-**based strategy", which utilizes the nearby pickup/drop-off locations to the customers' final origin/destination, has increasingly been studied and applied to improve the efficiency of DRT by reducing the operational costs with little inconvenience to customers (Czioska et al., 2019; Ma et al., 2021). Despite its benefit, the integration of an MP-based strategy into an electrified DRT system complicates vehicle dispatching and routing operations due to the interactions between customer-to-meeting point assignment, vehicle routing, and charging scheduling. Existing studies often assume unlimited charging station capacity, usually violated in practice because the availability of (fast-)charging infrastructure is very limited due to high investment costs. While electric vehicle routing problems have been studied extensively, efficient solution algorithms for solving medium/large problem instances considering capacity-constrained charging stations are still underdeveloped.

To fill those gaps, our previous study (Ma et al., 2023a,b) proposed a model for meeting-point-based electric feeder service with charging synchronization constraints (MP-EFCS), where the service radius of meeting points is a system parameter to trade off user's inconvenience (walking time) and operational costs (vehicle routing time and charging time). The developed method allows the system to jointly optimize the customer-to-meeting-point assignment and vehicle routing costs under vehicle charging synchronization constraints with partial recharges. The problem is a variant of dial-a-ride problems (Cordeau, 2006) with EVs, but is more complex due to the charging synchronization constraints. Figure 1 illustrates an example of the MP-EFCS service. Unlike traditional door-to-door DRT, in the MP-EFCS, customers are assigned to the one of predefined meeting points within a given maximum walking distance. Vehicles pick up customers at these meeting points and drop them off at their desired transit stations within desired arrival time windows.

Using a larger meeting point separation distance may reduce the number of customers due to the longer walking distance. In addition, when other system parameters are considered (i.e. fleet size, different weights used in the objective function, etc.), there might be synergies to optimize the system performance. However, these aspects have not been studied yet. Therefore, in this study, we provide new contributions to analyze the effect of different system parameters and configure the system parameters of MP-EFCS service in a study area. The main contributions are as follows.

- We analyze the impact of various system parameters of the MP-EFCS service, including fleet size, meeting point separation distance, and user-defined coefficients in the objective functions.
- Based on the identified key system parameters, we optimize their configuration to minimize the operator’s cost with the least impact on customer inconvenience.
- The computational experiments are performed for a morning peak-hour commuting scenario from Arlon to Luxembourg. The key performance indicators are analyzed to evaluate the performance of the MP-EFCS service. Its performance is compared with traditional dial-a-ride service.

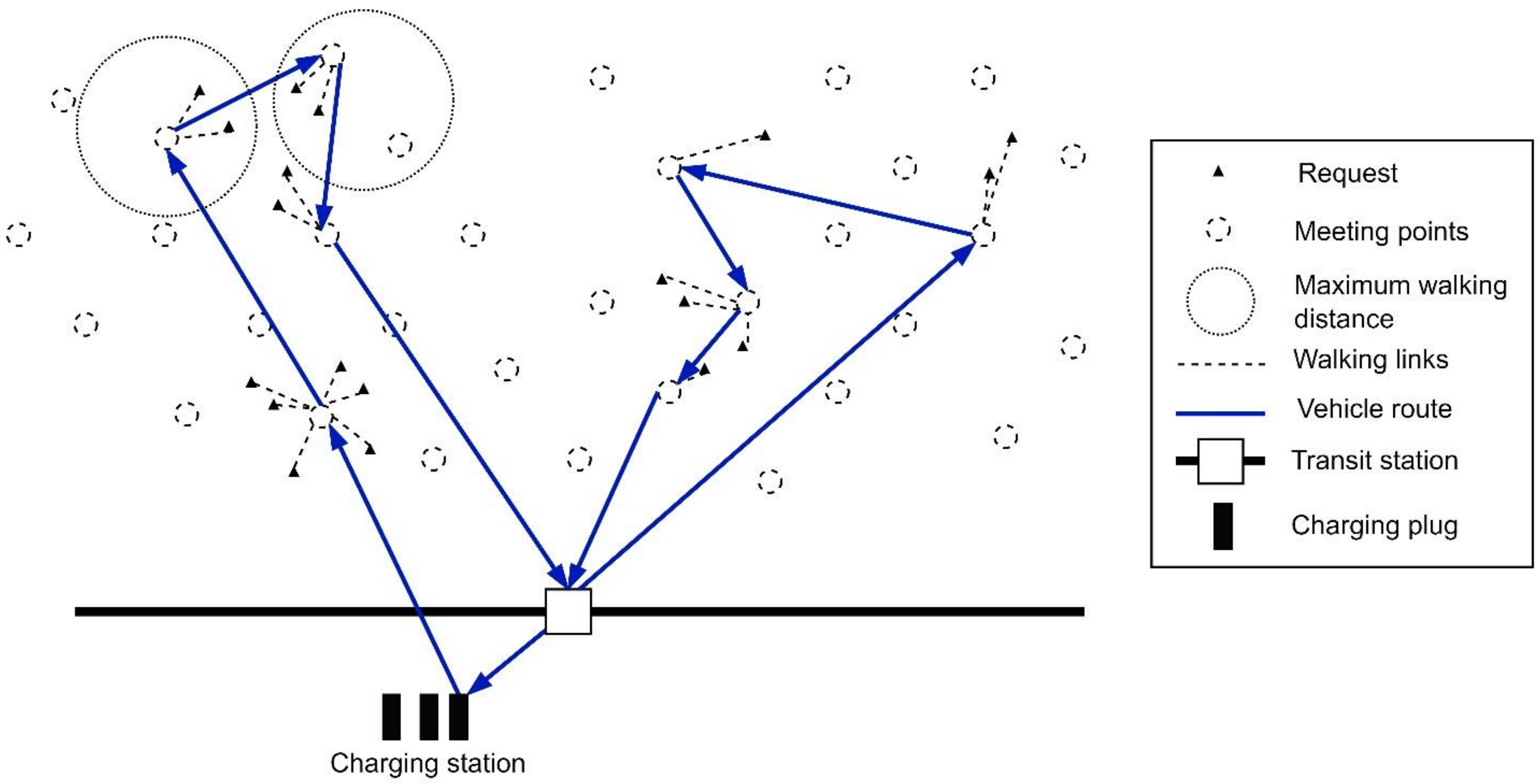

Figure 1. An illustrative example of an electric meeting-point-based feeder service. A vehicle makes two round trips from the transit station, picking up customers on the way and dropping them off at the transit station.

## 2. Related work

Demand-responsive transport (DRT) services have been studied and implemented since the seventies (Wilson et al., 1976). Extensive research has been conducted on this topic during the last half-century. The reader is referred to the recent works of Vansteenwegen et al. (2022) and Ma & Fang (2022) for a comprehensive review. In this section, we briefly review the meeting-point-based approach to improve the efficiency of classical door-to-door DRT and the methodology to configure the fleet size and the system parameters of DRT services.

The meeting-point-based DRT services consist of operating vehicle routes on a set of predefined “meeting points”, i.e., street corners, existing bus stops, or parking lots near the customer's origin or destination (Czioska et al., 2019). This approach has been studied for ridesharing services to allow for reducing vehicle routing costs to match the origins and destinations of drivers and riders (Stiglic et al., 2015). Similar concepts have been studied to reduce the operating costs of door-to-door DRT services. For example, Montenegro et al. (2022) propose a flexible feeder service in which a fleet of buses operates on a set of mandatory bus stops on a backbone line while allowing flexible bus routes to visit optional bus stops in the service area. The authors formulate the problem as a MILP and develop a large neighborhood search algorithm to solve it for large test instances (Montenegro et al., 2021). The analytical approach has been developed in the past for this type of flexible bus service (Chang and Schonfeld, 1991). A similar on-demand bus system is proposed by Melis and Sörensen (2022). The authors develop a MILP model to assign customer requests to bus stops in the first stage and then optimize bus routes in the second stage to minimize the total user ride time. An efficient large neighborhood search heuristic is developed to solve large test instances. The results show that the total user ride time is lower for on-demand bus systems compared to fixed-route bus services.

The fleet size and mix vehicle routing problem consists of determining the vehicle routing, type, and number of vehicles to serve a set of customer demands at minimum cost (Brandão, 2009). The problem can be formulated as a MILP with vehicle type-specific operating constraints. Since the problem is NP-hard, many heuristic and metaheuristic approaches have been developed (Hiermann et al., 2016). Hiermann et al. (2016) added different types of EVs to the problem, aiming to reduce both vehicle purchasing and routing costs. Other studies used an agent-based simulation approach for fleet size minimization considering stochastic demand (Scheltes & Correia, 2017; Winter et al., 2018). These studies focus on minimizing the total system cost by considering fleet utilization, maintenance, insurance, and/or depreciation costs. However, to achieve optimal system performance, other key system parameters need to be jointly configured. For electric vehicle routing problems with capacitated charging stations, the reader is referred to Lam et al. (2022) and Ma et al. (2023b) for related literature reviews.

Table 1 provides a summary and comparison of different DRT systems. The main distinguishing features of the MP-EFCS service in this study are that it is a more flexible meeting-point-based service and allows customers to be rejected. In this way, the system is more flexible and cost-efficient. For customer inconvenience, constraints related to maximum walking distance, maximum user riding time, and time window constraints at drop-off transit stations are considered. In addition, the system is operated with a fleet of EVs, and their charging operations are optimized with partial recharge policy and charging synchronization constraints.

Table 1. Characteristics of different demand-responsive transport systems

| Problem | DARP[1] | ODBRP[2] | DRFS[3] | MP-EFCS[4] |
|---|---|---|---|---|
| Objective fun. (Minimize) | Total routing cost | Total user ride time | The weighted sum of bus routing time, passenger walking time, and penalty for early/late arrivals at destinations | The weighted sum of bus routing time, charging time, passenger’s walking time, the penalty for early arrivals at destinations, and penalty for customer rejections |
| Max. ride time constr. | ✓ | ✓ | | ✓ |
| Time window constr. | Pickup or drop-off loc. | Drop-off loc. | Drop-off loc. | Drop-off loc. |
| Capacity constr. | ✓ | ✓ | ✓ | ✓ |
| Max. route duration constr. | ✓ | | | |
| Bus stop assignment | | ✓ | ✓ | ✓ |
| Max. walking dist. constr. | | | ✓ | ✓ |
| Customer rejection | | | | ✓ |
| EV | | | | ✓ |

| | |
|---|---|
| Charging capacity constr. | ✓ |
| Partial recharge | ✓ |

Remark: 1:DARP: dial-a-ride problem (Cordeau, 2006); 2. ODBRP: on-demand bus routing problem (Melis and Sörensen, 2022); 3. DRFS: demand-responsive feeder service (Montenegro et al., 2022); 4. This study.

## 3. Methodology

### 3.1. Problem description

Consider a meeting-point-based electric (first-mile) feeder service using a fleet of electric vehicles (or buses interchangeably) in a rural area. Customers are assumed to be willing to walk to meeting points within a maximum walking distance. Customers submit travel requests a day ahead via a dedicated platform detailing their origin, drop-off transit stations, and desired arrival time at transit stations which corresponds to the transit departure time. Based on the collected requests, the operator confirms requests, their pick-up time, and recommended meeting points or reject requests. Vehicles must reach the transit stations within a pre-defined buffer time (i.e., 10 minutes) before the departure of trains (transit service). Vehicles' charge levels have maximum and minimum bounds, with charging done at limited operator-owned chargers. The charging speed is assumed linear. Charging operations cannot overlap at any charger. The objective of the MP-EFCS service is to minimize the system costs. Four groups of constraints are considered as follows.

- **Customer-to-meeting-point assignment constraints**: Assign each customer to one meeting point within their maximum walking distance.
- **Vehicle routing constraints**: Ensure that the vehicle capacity is respected and that the vehicle arrival time and starting time of service at each node are consistent, given the associated time window and ride time constraints.
- **Energy constraints**: Specify the state of charge of vehicles at each node when traversing arcs and when performing charging operations.
- **Charging scheduling (synchronization) constraints**: Prevent simultaneous charging at any single charger

The MILP formulation of the MP-EFCS is as follows (Ma et al., 2023a).

$$Min\, Z = \lambda_1 \sum_{k \in K} \left( \sum_{(i,j) \in \mathcal{A}_B} t_{ij}\, x_{ij}^k + \sum_{s \in S'} \tau_s^k \right) + \lambda_2 \sum_{(r,i) \in \mathcal{A}_c} \sum_{k \in K} y_{ri}^k t_{ri} + \lambda_3 \sum_{k \in K} \sum_{i \in D'} W_i^k + \lambda_4 \sum_{(r,i) \in \mathcal{A}_c} \left( (1 - \sum_{k \in K} y_{ri}^k \right) \tag{1}$$

The objective function (1) seeks to minimize the weighted sum of vehicle travel time ($t_{ij}$) and charging time ($\tau_s^k$) (first term), customer walking time ($t_{ri}$) (second term), excess vehicle waiting time at transit stations ($W_i^k$) (third term), and the penalty for unserved customers (fourth term). $\lambda_1$ to $\lambda_3$ are user-defined weights to account for the tradeoffs between operational costs and customer incontinence. $\lambda_4$ is the penalty for rejecting customers. The notation is shown in Appendix A.

Customer-to-meeting-point assignment constraints:

$$\sum_{k \in K} \sum_{i \in G'} y_{ri}^k \leq 1, \qquad \forall r \in R \tag{2}$$

$$\sum_{k\in K}\sum_{i\in G'} w_{ri}y_{ri}^k \le w_{max}, \qquad \forall r \in R \tag{3}$$

Vehicle routing constraints:

$$\sum_{j\in G'\cup S'\cup\{N+1\}} x_{0j}^k = 1, \forall k \in K \tag{4}$$

$$\sum_{i\in\{0\}\cup S'\cup D'} x_{i,N+1}^k = 1, \forall k \in K \tag{5}$$

$$\sum_{j\in G'\cup\{N+1\}} x_{sj}^k \le 1, \qquad \forall k \in K, s \in S' \tag{6}$$

$$\sum_{i\in V_0} x_{ij}^k \le 1, \qquad \forall k \in K, j \in G' \tag{7}$$

$$\sum_{i\in V_0} x_{ij}^k - \sum_{i\in V_{N+1}} x_{ji}^k = 0, \qquad \forall k \in K, j \in V \tag{8}$$

$$\sum_{r\in R} y_{ri}^k \le M_1 \sum_{j\in V_{\mathrm{N+1}}} x_{ij}^k, \qquad \forall k \in K, i \in G' \tag{9}$$

$$y_{ri}^k = 1 \Rightarrow \sum_{j\in V_0} x_{ji}^k = \sum_{j\in G'\cup D'} x_{jd_r}^k, \qquad \forall k \in K, i \in G', r \in R \tag{10}$$

$$x_{ij}^k = 1 \Rightarrow q_j^k = q_i^k + \sum_{r\in R} y_{rj}^k, \qquad \forall k \in K, i \in V_0, j \in G' \tag{11}$$

$$x_{ij}^k = 1 \Rightarrow q_j^k = q_i^k - \sum_{r\in R}\sum_{g\in G'} y_{rg}^k, \qquad \forall k \in K, i \in G', j \in D' \tag{12}$$

$$0 \le q_i^k \le Q^k, \qquad \forall k \in K, i \in V_{0,N+1} \tag{13}$$

$$B_j^k \ge B_i^k + u_i + t_{ij} - M_2(1 - x_{ij}^k), \qquad \forall k \in K, i \in V_0, j \in V_{N+1} \tag{14}$$

$$B_j^k \ge B_s^k + \tau_s^k + t_{sj} - M_2(1 - x_{sj}^k), \qquad \forall k \in K, s \in S', j \in \{G' \cup N+1\} \tag{15}$$

$$x_{ij}^k = 1 \Rightarrow A_j^k = B_i^k + t_{ij} + u_i, \qquad \forall k \in K, i \in G' \cup D', j \in D' \tag{16}$$

$$W_i^k \ge B_i^k - A_i^k - M_2(1 - p_i^k), \qquad \forall k \in K, i \in D' \tag{17}$$

$$p_i^k = \sum_{j\in V} x_{ji}^k, \qquad \forall k \in K, i \in D' \tag{18}$$

$$A_{d_r}^k - B_i^k - u_i \le L_i + M_2\left(1 - y_{ri}^k\right), \qquad \forall k \in K, (r,i) \in A_C \tag{19}$$

$$e_i \le B_i^k \le l_i, \qquad \forall k \in K, i \in D' \tag{20}$$

Vehicle energy constraints:

$$E_0^k = E_{init}^k, \qquad \forall k \in K \tag{21}$$

$$E_{min}^k \le E_i^k \le E_{max}^k, \qquad \forall k \in K, i \in V \tag{22}$$

$$x_{ij}^k = 1 \Rightarrow E_j^k = E_i^k - \beta^k c_{ij}, \qquad \forall k \in K, i \in V_0 \backslash S', j \in V_{N+1} \tag{23}$$

$$x_{ij}^k = 1 \Rightarrow E_j^k = E_s^k + \alpha_s \tau_s^k - \beta^k c_{sj}, \qquad \forall k \in K, s \in S', j \in \{G' \cup N+1\} \tag{24}$$

Charging scheduling constraints:

$$v_s = \sum_{k\in K}\sum_{j\in G'\cup N+1} x_{sj}^k, s \in S' \tag{25}$$

$$v_h \le v_l, \forall h, l \in S_o', o \in S, h < l \tag{26}$$

$$\sum_{k\in K} B_h^k \ge \sum_{k\in K} B_l^k + \sum_{k\in K} \tau_l^k - M_2(2 - v_h - v_l), \forall h, l \in S_o', o \in S, h < l \tag{27}$$

$$\tau_s^k + B_s^k \leq M_2 \sum_{j \in G' \cup N+1} x_{sj}^k, \forall s \in S', k \in K \tag{28}$$

$$v_s \leq 1, s \in S' \tag{29}$$

Domain of variables:

$$x_{ij}^k \in \{0,1\}, \qquad \forall k \in K, i,j \in V_{0,N+1} \tag{30}$$

$$y_{ri}^k \in \{0,1\}, \qquad \forall k \in K, r \in R, i \in G' \tag{31}$$

$$\tau_s^k \geq 0, v_s \geq 0, \forall k \in K, s \in S' \tag{32}$$

$$A_i^k \geq 0, B_i^k \geq 0, \forall k \in K, i \in V_{0,N+1} \tag{33}$$

$$p_i^k \in \{0,1\}, W_i^k \geq 0, \forall k \in K, i \in D' \tag{34}$$

The operator ⇒ means that if the left-hand side of the equations is true, then their right-hand side constraints need to be satisfied. Equation (2) states that each customer can be assigned to at most one meeting point. Equation (3) imposes a customer's maximum walking distance to meeting points. Equations (4)-(5) and (8) are the flow conservation. Equations (6)-(7) state that each meeting point and each charger can be visited by a vehicle at most once. Equation (9) ensures that a meeting point assigned to a customer must be visited by a vehicle. Equation (10) ensures that if customer $r$ is assigned to the meeting point $i$, the drop-off location of customer $r$ (i.e. $d_r$) must be served by the same vehicle $k$. Equations (11)-(13) manage the passenger loads at meeting points and transit stations under vehicle capacity. Equations (14)-(15) ensure the consistency of the beginning of service time at bus nodes and charger nodes with a charging duration of $\tau_s^k$, respectively. Equation (16) defines the relationship between vehicle arrival time and the beginning time of service. Equations (17)-(18) determine the vehicle waiting time at transit stations when arriving earlier than its buffer time. Equation (19) limits the customer's maximum riding time (i.e. 1.5 times the direct riding time). Equation (20) sets up the arriving time constraints (within a fixed buffer time before the departure time of trains) of vehicles at transit nodes. Equations (21)-(22) initiate the state of charge of vehicles and their upper and lower bounds. Equations (23)-(24) are energy conservation constraints. Equations (25)-(29) are charging scheduling constraints with capacitated charging stations. Equations (30)-(34) specify the domain of variables.

Since the problem is NP-hard, exact solutions are possible only for small instances. Thus, we previously developed a hybrid deterministic annealing (DA) metaheuristic for larger instances (Ma et al. 2023a, b). The proposed hybrid metaheuristic consists of three steps: 1. assign customers to nearby meeting points (see Appendix B); 2. optimize electric vehicle routing and charging scheduling with partial recharge policy and charging synchronization; 3. re-optimize the solution obtained in Step 2 to insert unserved customers and re-optimize current routes based on a metaheuristic approach. The algorithm was tested on 20 test instances with up to 100 customers and 49 used meeting points under different initial battery levels and demand distributions (peaked and non-peaked demand scenarios). The proposed algorithm efficiently found good-quality solutions in a short computational time. Compared to a 4-hour solution from the MILP solver, the average and best solution gaps are 4.02% and 3.31% respectively. The reader is referred to the aforementioned studies for a more detailed description. In the following section, we apply the developed metaheuristic to a real-world case study. The focus is on the effect of different system parameters and their configuration to obtain a better performance of the system.

### 3.2. Case study design

A case study mimics a real-world scenario for Belgian cross-border workers commuting from Arlon to Luxembourg during the morning peak hour (see Figure 2). We assume that the MP-EFCS service is provided in the Arlon region to connect customers to the main train station with one railroad line (line 50) to Luxembourg. During morning peak hours (6:00-9:00), trains depart every 15 minutes on average. **Five** test instances are generated in this context with customer requests randomly distributed within a radius between 1.5 km and 6km. The total number of requests is assumed as 600 from 6:00-9:00 which corresponds to ~50%[1] of total train trips for the Arlon to Luxembourg cross-border workers based on the 2017 Luxmobil survey (Lambotte et al., 2021). Customers' desired arrival times at the station corresponding to the train departure time follow a normal distribution and are randomly generated for each test instance. With a maximum walking distance of customers of 1.0 km, no customers are assumed to walk directly to the station. Using a grid system, meeting points are set 1.2 km apart for the base scenario. The characteristics of the MP-EFCS service are shown in Table 2. To facilitate the analysis, a homogenous 24-seat shuttle is used. However, the proposed method can be adapted by considering a heterogeneous fleet with multiple depots without difficulty. We assume that the train station and the depot are located at (0, 0) and the two DC fast chargers are located at (0, 0) and one at (1, 0) (in km), respectively (see Figure 3). Note that the impact of different charging infrastructure configurations will be analyzed in our future work. To invoke charging operations, initial states of charge of vehicles are set randomly between 50% and 100% of battery capacity.

Table 2. Summary of the parameter setting of the case study.

| Parameter | Value | Parameter | Value |
|---|---|---|---|
| Passenger Capacity | 24 | Maximum walking distance (km) | 1.0 |
| Battery Capacity (kWh) | 118 | Nearest separation distance between meeting points (km) | 1.2 |
| Energy consumption rate (kWh/km) | 1.23 | Walking speed (km/hour) | 5.1 |
| Bus speed (km/hour) | 30 | Detour factor | 1.5 |
| Number of chargers | 3 | Charging power (kWh/hour) | 50 |

*Based on the characteristics of the EVSTAR shuttle (https://greenpowermotor.com/wp-content/uploads/Brochures/EVSTAR_Brochure.pdf)

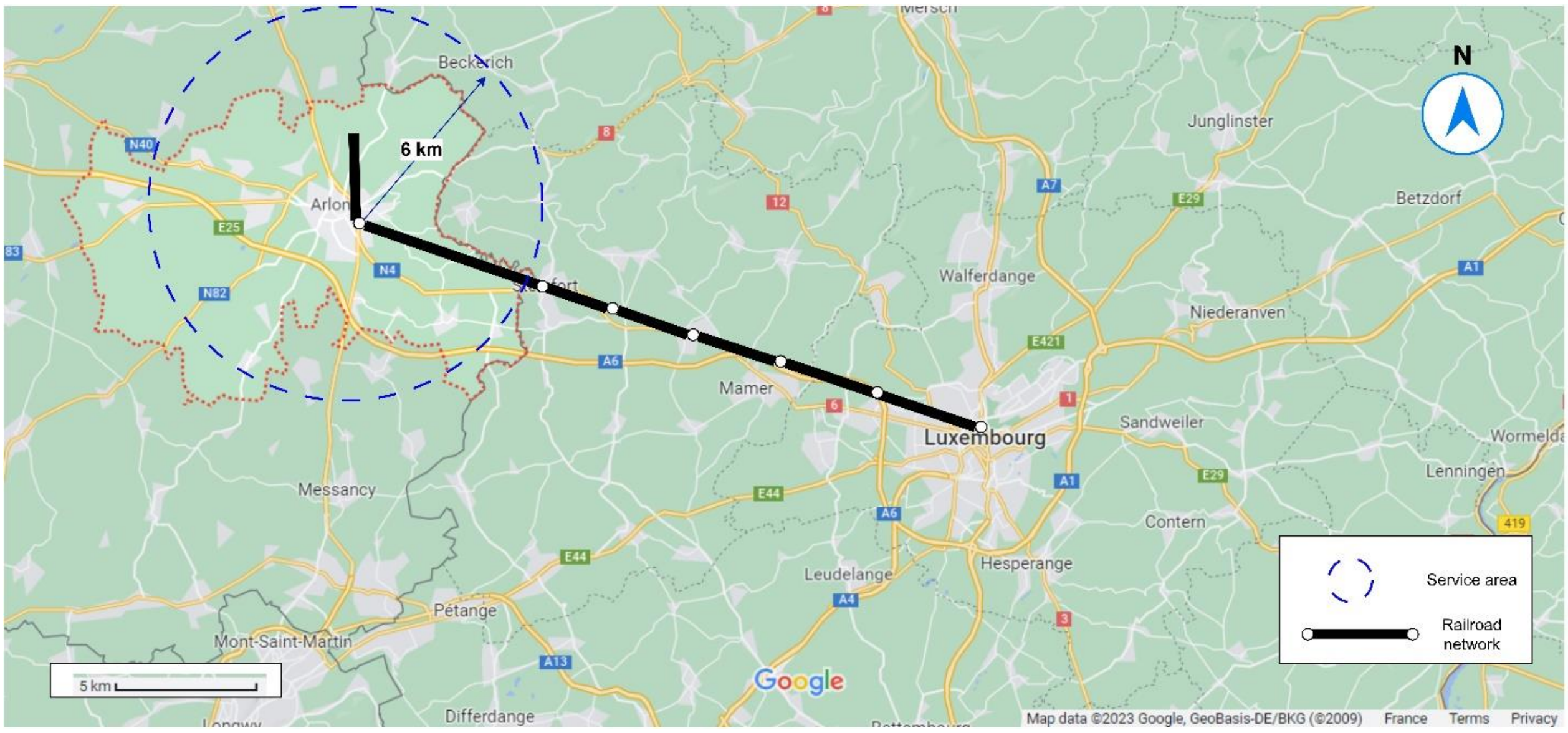

Figure 2. Feeder service area in the region of Arlon and its railway connection to Luxembourg.

[1] Based on the number of Arlon's cross border workers (16900 i.e. 9% of Luxembourg's cross border workers in 2023) multiplied by the mode share of train for this population (7.9%). As a result, a total of 1335 train trips on weekdays are estimated for the cross border workers from Arlon to Luxembourg.

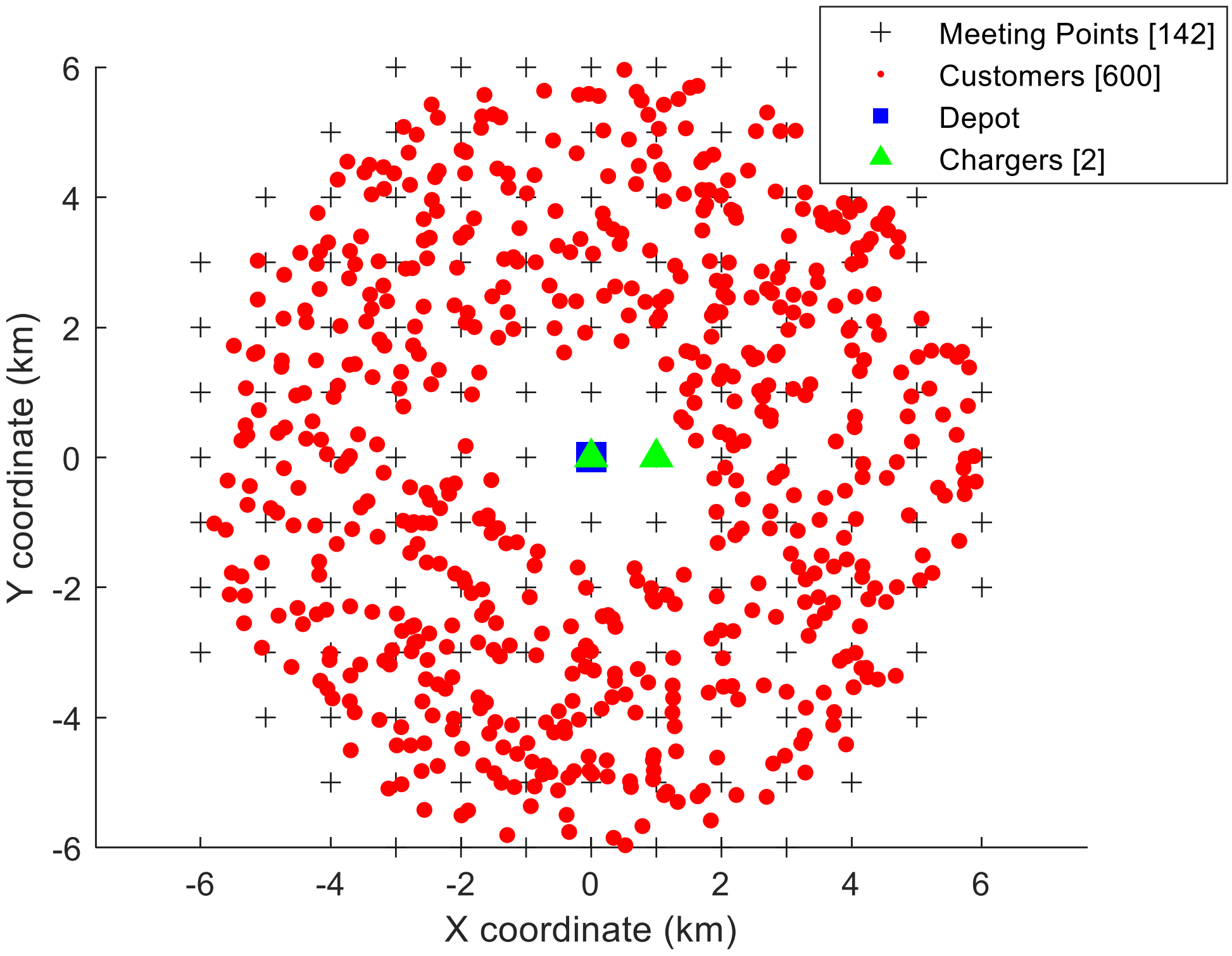

Figure 3. Example of a test instance with 600 requests during the morning peak hour in the study area (meeting point separation distance is 1.0 km).

### 3.3. Computational experiment settings

We design four computational experiments to investigate the effects of different parameters on the performance of the system. A set of key performance indicators (KPIs) are used to evaluate the system performance, including the number of used vehicles, service rate, average walking distance, average in-vehicle travel time (IVT), additional vehicle waiting time at transit stations (earlier arrivals before the buffer time), total kilometers traveled by vehicles (KMT), the average number of served customers per KMT, average customers per meeting point, and total charging time of the fleet.

These experiments first evaluate the effect of fleet size, meeting point availability, and user-defined weights (in the objective function) on KPIs, while keeping other system parameters unchanged based on the base case (described below). In doing so, we can isolate the effect when changing the value of one single system parameter. Then we configure these parameters over a set of values to improve the performance of the system compared to the base case. All experimental results are based on the average performance of the five test instances. For each test instance, a 3-run best solution is used. Following our previous study (Ma et al., 2023b), the same tuned algorithmic parameters of the metaheuristic are used for this study.
The set of tested values is set as follows.

- Fleet size: $|K| \in \{10,12,14,16\}$
- Meeting point separation distance (in meters): $\tilde{d} \in \{1000, 1200, 1400\}$
- Weights in the objective function: $\lambda_1 \in \{0, 1, 2, 3\}, \lambda_2 \in \{0, 1,2,3\}, \lambda_3 \in \{0,1,2,3\}$

The penalty for rejecting customers ($\lambda_4$) is set as 40 for all experiments. This is based on our preliminary studies showing that using a higher penalty results in similar performance. The parameters of the **base case** are set as $\lambda_1 = \lambda_2 = \lambda_3 = 1$, the meeting point separation distance = **1.2 km**, and the fleet size = **14 vehicles**. Note that the effect of other system parameters can be analyzed similarly as well (e.g. buffer time at transit stations, detour factor for specifying riding time constraints, etc.).

The computational experiments are conducted on a PC with Intel(R) Core(TM) i7-11800H, 16 logical processors, and 64 GB of memory. The implementation of the metaheuristic algorithm is based on Julia and uses Gurobi MIP solver version 10.0.2.

## 4. Results

### 4.1. The effect of fleet size

The effect of fleet size on the KPIs is shown in Table 3 and Figure 4. The minimum fleet size to serve more than 95% of requests is 14 vehicles, given the base case parameter setting. When reducing the fleet size from 14 to 12 and 10, the service rate (i.e. percentage of the number of served customers) drops from 96.5% to 89.2% (-7.3%) and 78.2% (-18.3%), respectively and total charging time (in minutes) of the fleet increases from 38.1 minutes to 64.1 and 71.6 minutes, respectively (see Table 3). This is because when more vehicles are available, the average kilometer traveled per vehicle decreases so the charging needs tend to decrease with shorter routing distance per vehicle. As the cost to use an additional fleet is not considered, all available vehicles are used for all cases. Average IVT and walking distance do not vary significantly over different fleet sizes. The average number of customers per meeting point over different fleet sizes is 2.83, while number of served customers per KMT is between 0.65 (16 vehicles) and 0.75 (10 vehicles).

Table 3. KPIs of the MP-EFCS with different fleet sizes.

| Fleet size | # used vehicles | Service rate (%) | Walking dist (km) | IVT | Charging time | KMT | Cus /KMT | Cus /MP | CPU (sec) |
|---|---|---|---|---|---|---|---|---|---|
| 10 | 10 | 78.2 | 0.58 | 9.68 | 71.6 | 626.8 | 0.75 | 2.82 | 1317 |
| 12 | 12 | 89.2 | 0.57 | 9.83 | 64.1 | 752.3 | 0.71 | 2.79 | 1252 |
| 14 | 14 | 96.5 | 0.58 | 9.90 | 38.1 | 859.8 | 0.67 | 2.81 | 988 |
| 16 | 16 | 99.6 | 0.59 | 10.03 | 13.3 | 915.6 | 0.65 | 2.88 | 771 |

Remark: IVT: in-vehicle travel time (minute); KMT: Total kilometers traveled by vehicles (km); cus/KMT: average number of served customers per KMT; cus/MP: average number of customers of activated meeting points. Excess waiting time at transit stations is 0 for all tested cases.

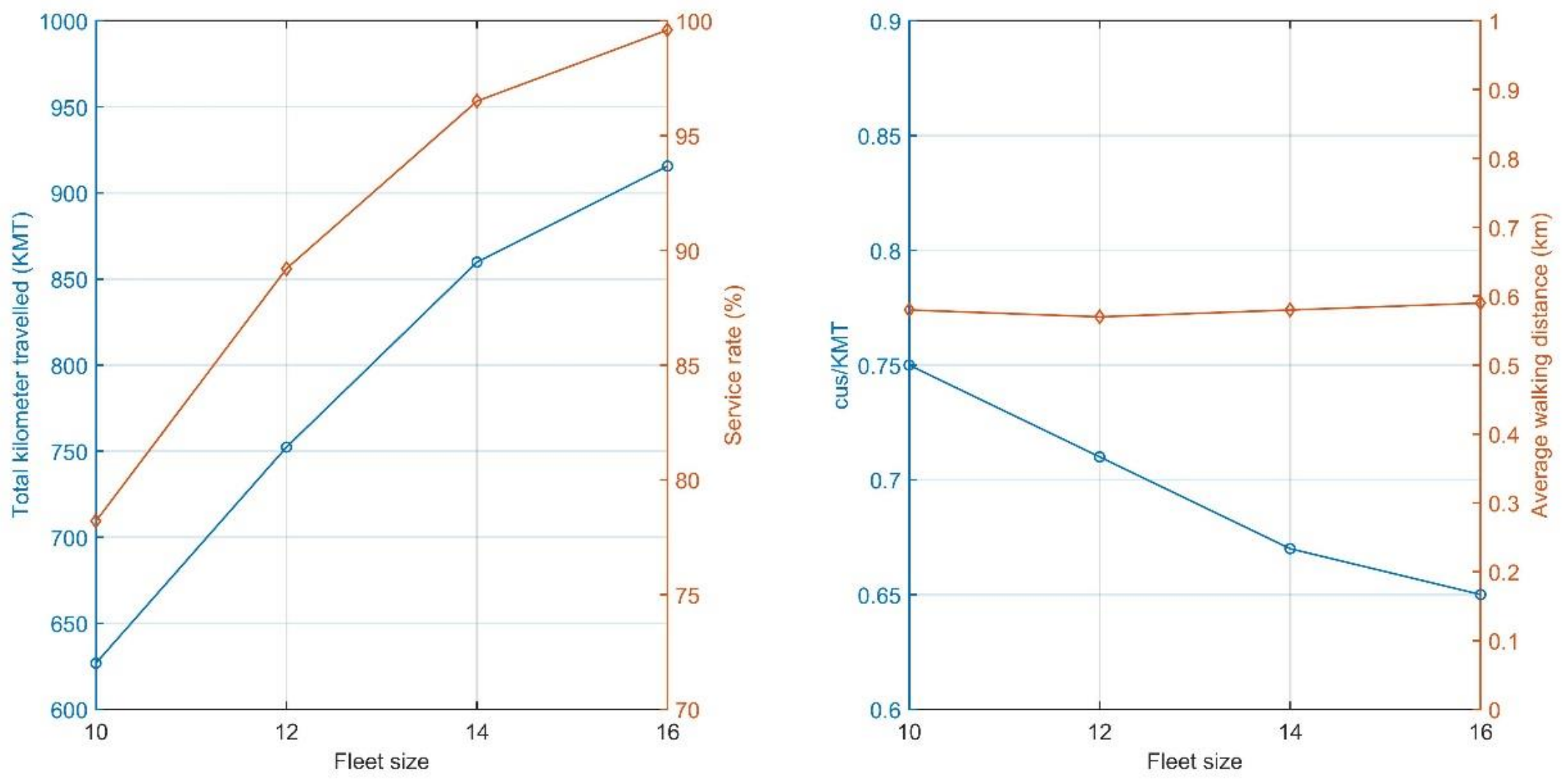

Figure 4. The comparison of KPIs over different fleet sizes.

### 4.2. The impact of meeting point availability

We change the MP separation distance from 1000 m to 1400 m, giving the base-case parameter setting. Note that in our experiments, the meeting points (MPs) are organized as a grid system (see the example in Figure 3). An $X$km MP separation distance results in the furthest away from the grid center of $X/\sqrt{2}$ km. So the maximum walking distance of 1 km will not be an issue for the 1400 m case ($1.4/\sqrt{2} = 0.9899$km). However, the MP separation distance will change the number of MPs within walking distance and hence the flexibility of pickup locations. As a result, more customers are served when there are more MPs (e.g., 99.6% are served with MP dist = 1000m) than the case with fewer MPs (e.g., 96.3% are served when MP dist = 1200m).

The smaller MP separation distances (i.e., more MPs) result in shorter KMTs (-1.62% and -5.44% compared to using 1200 m and 1400 m, respectively). It is because as more MPs become available, the chances for multiple customers to share the same MP increases, and the average number of customers per MP (cus/MP) increases (see Table 4). Hence, the number of used MPs decreases leading to less KMT as well as less total charging time (see Table 4). These results indicate that more available MPs allow the system to reduce operational costs by jointly optimizing customer-to-meeting-point assignment and vehicle routing.

Table 4. Key metrics of the MP-EFCS with different separation distances between meeting points.

| MP dist | # of MP | Service rate (%) | Walking dist | IVT | Charging time | KMT | cus/KMT | cus/MP | CPU (sec) |
|---|---|---|---|---|---|---|---|---|---|
| 1000 | 141 | 99.6 | 0.59 | 9.95 | 9.41 | 822.89 | 0.73 | 3.14 | 1093 |
| 1200 | 100 | 96.3 | 0.58 | 9.92 | 30.62 | 856.10 | 0.68 | 2.81 | 1331 |
| 1400 | 74 | 91.67 | 0.57 | 10.03 | 41.13 | 870.24 | 0.63 | 2.56 | 949 |

Remark: MP dist: Distance between two closest adjacent meeting points. The excess waiting time is 0 for all tested cases.

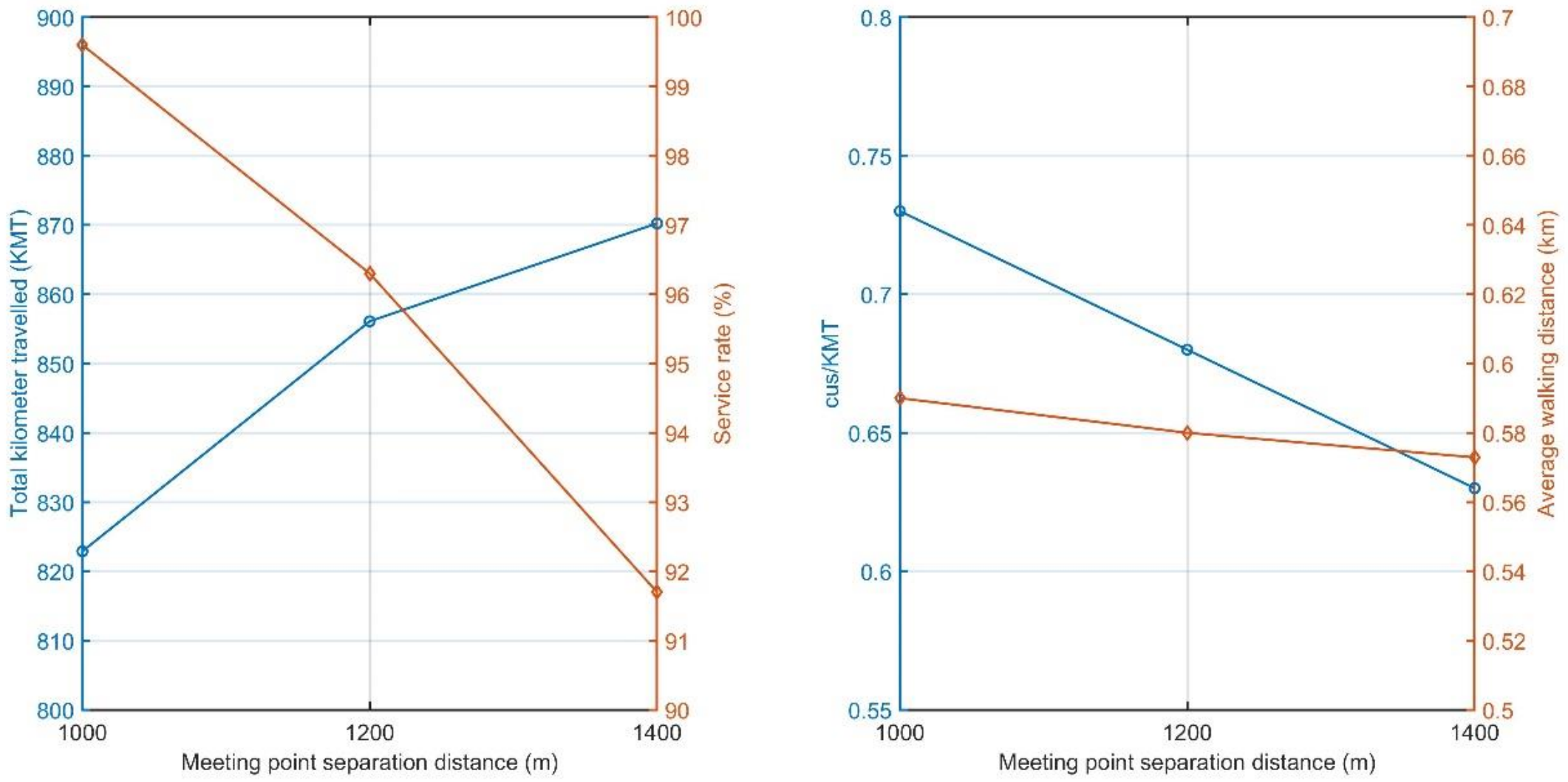

Figure 5. The comparison of KPIs over different meeting point separation distances.

### 4.3. The effects of the objective function parameters

We examine various user-defined weights in the objective function: $\lambda_1$ is related to total vehicle routing and charging times, while $\lambda_2$ and $\lambda_3$ addresse walking distance and excess waiting time at a transit station. Parameter set 1 (see Table 5) only considers the total routing time and charging time; set 2 minimizes total walking time and excess waiting time at transit stations; set 3 is the base case with equal weights. Parameter sets 4 and 5 increase the weights of $\lambda_1$ in favor of operators. Table 5 shows the effect of $\lambda_1$ on the system performance. First, solely minimizing customers' walking distance and excess waiting time of vehicles at transit stations results in the highest KMT with significant unserved customers (serve rate is only 71.6%, the line of parameter set 1 on Table 5). Conversely, only minimizing the total vehicle routing time and charging time results in the highest average walking distance (0.61 km). Compared to the base case, parameter set 5 gives the lowest KMT (813.9 km) and significantly lower total charging time (2.8 minutes), suggesting both the operator's and customer's costs need to be considered to have synergies.

Table 5. Effects of $\lambda_1$.

| Para. set | $\lambda_1$ | $\lambda_2$ | $\lambda_3$ | Service rate (%) | Walking dist | IVT | Charging time | KMT | cus/KMT | cus/MP | CPU (s) |
|---|---|---|---|---|---|---|---|---|---|---|---|
| 1 | 0 | 1 | 1 | 71.6 | 0.46 | 9.87 | 54.1 | 877.0 | 0.49 | 1.68 | 1397 |
| 2 | 1 | 0 | 0 | 97.3 | 0.61 | 9.91 | 26.7 | 849.6 | 0.69 | 2.80 | 999 |
| 3 | 1 | 1 | 1 | 96.5 | 0.58 | 9.90 | 38.1 | 859.8 | 0.67 | 2.81 | 988 |
| 4 | 2 | 1 | 1 | 97.9 | 0.58 | 9.97 | 10.6 | 834.0 | 0.70 | 2.89 | 872 |
| **5** | **3** | **1** | **1** | **97.2** | **0.58** | **9.97** | **2.8** | **813.9** | **0.72** | **2.85** | **1013** |

Remark: Fleet size is 14 and $\lambda_4$=40. Excess waiting time at the transit station is 0 for all tested cases.

For the effects of $\lambda_2$ (Table 6), the result is as expected that increasing $\lambda_2$ results in shorter average walking distance of customers. Regarding $\lambda_3$, there are no significant effects on the system performance compared with the base case.

Table 6. The effect of $\lambda_2$ and $\lambda_3$.

| $\lambda_1$ | $\lambda_2$ | $\lambda_3$ | Service rate (%) | Walking dist | IVT | Charging time | KMT | cus/KMT | cus/MP | CPU (s) |
|---|---|---|---|---|---|---|---|---|---|---|
| 1 | **0** | 1 | 97.6 | 0.61 | 9.78 | 47.0 | 866.3 | 0.68 | 2.83 | 843 |
| 1 | **1** | 1 | 96.5 | 0.58 | 9.90 | 38.1 | 859.8 | 0.67 | 2.81 | 988 |
| 1 | **2** | 1 | 97.1 | 0.58 | 9.98 | 26.1 | 847.9 | 0.69 | 2.85 | 859 |
| 1 | **3** | 1 | 95.3 | 0.57 | 9.87 | 30.2 | 853.9 | 0.67 | 2.81 | 992 |
| 1 | 1 | **0** | 96.5 | 0.58 | 9.97 | 36.1 | 860.3 | 0.67 | 2.77 | 1102 |
| 1 | 1 | **1** | 96.5 | 0.58 | 9.90 | 38.1 | 859.8 | 0.67 | 2.81 | 988 |
| 1 | 1 | **2** | 95.2 | 0.58 | 9.86 | 40.5 | 863.0 | 0.66 | 2.79 | 996 |
| 1 | 1 | **3** | 98.2 | 0.58 | 9.98 | 35.0 | 856.2 | 0.69 | 2.84 | 916 |

Remark: PS: Parameter set; Fleet size is 14 and $\lambda_4$=40. Excess waiting time at the transit station is 0 for all tested cases.

## 4.4. Optimizing system parameter configuration

Based on the above analysis, we search optimal parameter configuration concerning fleet size, meeting point separation distance, and $\lambda_1$ and analyze the trade-off of operational costs and customers' inconvenience. We considered:

- Fleet size: $|K| \in \{12,13,14\}$
- Meeting point separation distance: $\tilde{d} \in \{900, 1000, 1100\}$
- $\lambda_1$: $\lambda_1 \in \{1, 2, 3\}$(with $\lambda_2 = \lambda_3 = 1$)

Keeping other parameters as default, this gives 27 combinations. Table 7 summarises the best system parameter configuration only. Compared to the base case, using a smaller MP distance and higher $\lambda_1$ results in significant total KMT savings (-10.8% compared with the base case for 14 vehicles) with little increase in average walking time (0.62 or 0.63 km compared with 0.58 km for the base case). When reducing fleet size from 14 to 13 and 12, the configured system parameters allow for reducing further operational costs (-11.4% and -16.3% in terms of KMT) with +0.3% and -3.1% unserved customers. The result suggests that the best performance is achieved using a smaller MP distance with the joint configuration of system parameters to minimize the overall system costs.

Table 7. Best system parameter configurations for different fleet sizes.

| Best configuration | MP dist | # of used veh | $\lambda_1$ | Service rate | Walking dist | IVT | Charging time | KMT | cus/KMT | cus/MP |
|---|---|---|---|---|---|---|---|---|---|---|
| V14 | 900 | 14 | 3 | 100.0 | 0.63 | 9.92 | 0.0 | 763.6 (-10.8%) | 0.79 | 3.36 |
| V13 | 900 | 13 | 3 | 99.4 | 0.62 | 9.88 | 5.4 | 758.5 (-11.4%) | 0.79 | 3.35 |
| V12 | 900 | 12 | 3 | 96.0 | 0.63 | 9.89 | 27.5 | 716.6 (-16.3%) | 0.80 | 3.32 |
| Base case | 1200 | 14 | 1 | 96.3 | 0.58 | 9.92 | 30.6 | 856.1 | 0.68 | 2.81 |

Remark: Excess waiting time at the transit station is zero for all test cases. Based on the average performance of the 5 test instances with the 3-run best solution.

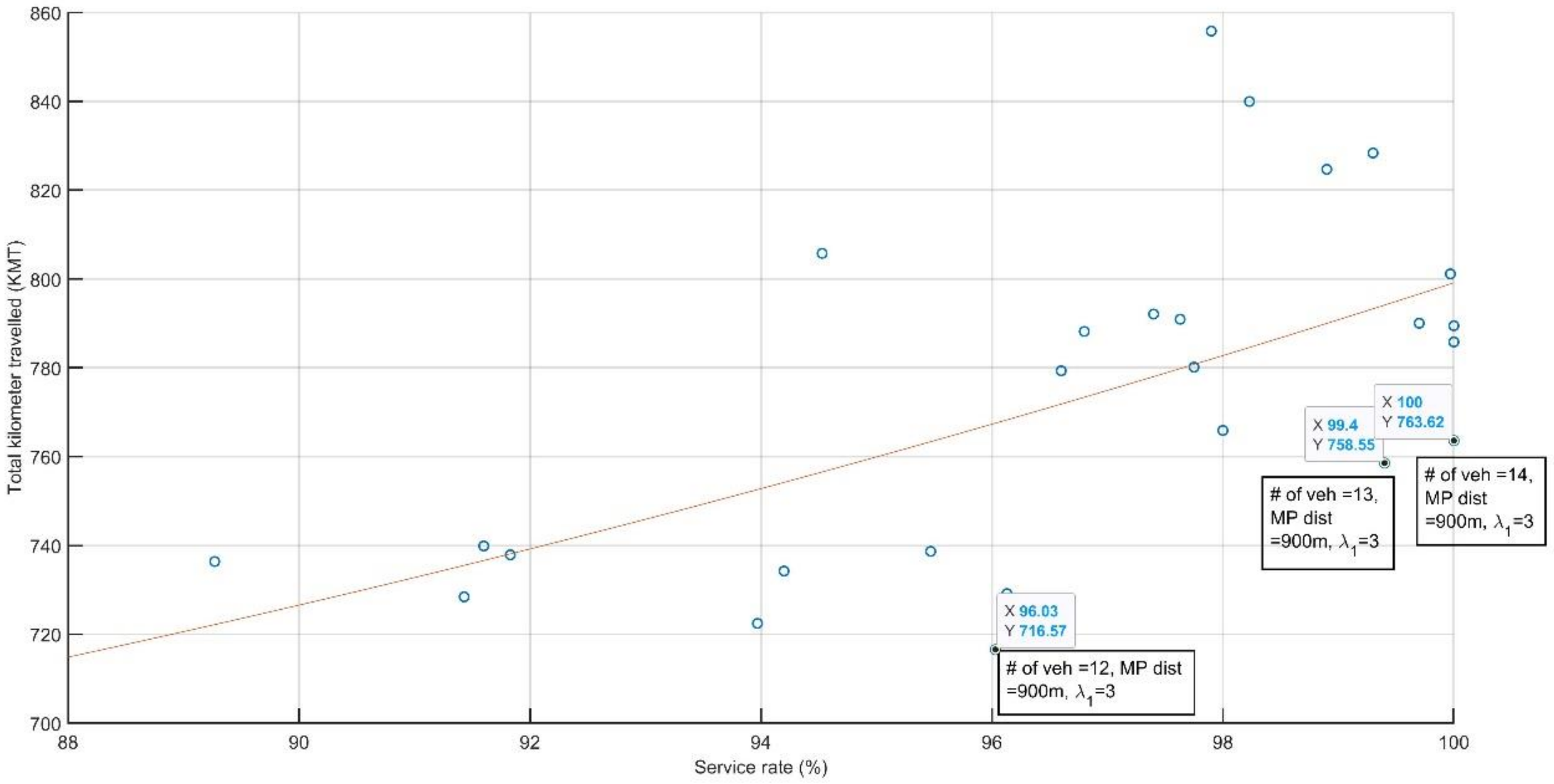

Figure 6. The distribution of the total KMT and service rates with respect to different System parameter configurations.

## 4.5. Benefits of the MP-EFCS service compared to door-to-door DARP

We further compare the performance of the MP-EFCS service with the door-to-door DARP. The objective function of the DARP minimize total travel time and charging time. Table 8 shows that while DARP using diesel vehicles (lower bound) requires 47 vehicles on average, MP-EFCS needs only 14 vehicles, a 70.2% reduction. Besides, the MP-EFCS service reduces 49.4% of KMT with an average of 0.79 customers per KMT, almost doubling the efficiency. Note that while the fleet use cost is not incorporated in the objective function, our metaheuristic attempts to reduce the number of used vehicles in the large neighborhood search process.

By converting the KMT saved by EVs into CO2 savings following Rosero et al.(2020), we estimate a yearly reduction of nearly 298 tons of CO2 emissions assuming that a fleet traveling 1000 km per day for 365 days (i.e., 365*1000*0.816kg = 297.84).

Table 8. Comparison of the KPIs between the MP-EFCS and door-to-door DARP

| Type of service | Service rate (%) | # of used vehicles | IVT | KMT | cus/KMT | cus/MP |
|---|---|---|---|---|---|---|
| DARP[1] | 100.0 | 47 | 10.6 | 1508.9 | 0.40 | - |
| MP-EFCS | 100.0 | 14 (-70.2%) | 9.92 (-6.4%) | 763.6 (-49.4%) | 0.79 (+97.5%) | 3.36 |

Remark: 1. Using diesel vehicles to have lower bounds.

## 5. Discussion and Conclusions

DRT often aims to improve fixed-route public transport in rural areas. However, DRT's high operational costs often lead to service terminations (Haglund et al., 2019). The push towards zero-emission transit and EV deployment further complicates DRT by demanding the joint optimization of routes and charging. Addressing these, our previous studies proposed a meeting-point-based EV-DRT model for a cleaner and more efficient service (Ma et al., 2023 a, b). In this study, we investigate the effects of different system parameters on a set of key performance indicators using a set of real-world size test instances. Simulation experiments revealed several findings and limitations of the proposed methodology.

a. The influence of meeting point availability is particularly important to reduce the number of unserved customers and routing costs (in terms of KMT). When the availability of meeting points is higher, the system allows to better configure customer assignment and routing jointly to enhance the system efficiency. While we conducted exploratory experiments, the operators could consider all possible meeting points (street corners, parking places, existing bus stops, etc.) in the service area to have the best performance of the system. When extending the problem to dynamic cases (i.e. allowing new requests to be inserted in the current routes after starting the service), dynamic meeting point assignment methodology can be applied (Czioska et al., 2019).
b. Evaluating fleet size, meeting point availability, and objective function weights simultaneously showed that an optimal configuration significantly reduces vehicles used and total KMT without compromising customer convenience. This suggests that this effort is necessary to optimize the system's performance and reduce operational costs.
c. When comparing the MP-EFCS services with traditional DARP, there is a clear advantage for this kind of meeting-point-based service. Our case study shows that it can reduce 70.2% of the number of used vehicles and 49.4% of KMT when the ride time constraint of customers is characterized by a detour factor of 1.5. This suggests a good potential for applying the MP-EFCS service to reduce the operational costs of traditional DRT services.
d. In terms of methodology, we show that the developed hybrid metaheuristic can solve real-world problems within reasonable computational time. To achieve better performance, we are currently improving the algorithm by allowing different customer-to-meeting-point assignment outcomes to be explored jointly and incorporating new local destroy-repair local search operators in the large neighborhood search process.

Future extensions include studying urban morphology's effect on the performance of the meeting-point-based feeder services; joint charging infrastructure and fleet size planning by considering both first- and last-mile feeder service with heterogeneous and mixed fleets (gasoline and EVs).

**Availability of data and material**: https://github.com/tym2021/datasets-for-electrified-meeting-point-based-feeder-bus-services.git

**Acknowledgements**: The work was supported by the Luxembourg National Research Fund (C20/SC/14703944).

**Appendix A. Notation.**

| | |
|---|---|
| *Sets* | |
| $G$ | Set of physical meeting points, i.e. $G = \{1, \ldots, N_G\}$ |
| $G'$ | Set of dummy (duplicate) meeting point vertices (nodes) |
| $D$ | Set of physical transit stations, i.e. $D = \{1, \ldots, N_D\}$ |
| $D'$ | Set of dummy (duplicate) transit station vertices |
| $S$ | Set of physical chargers, i.e. $S = \{1, \ldots, N_S\}$ |
| $S'$ | Set of dummy (duplicate) charger vertices |
| $R$ | Set of customers (i.e. location of origin of customers) |
| $K$ | Set of electric buses |
| $\bar{V}$ | Set of all vertices, i.e. $\bar{V} = G' \cup D' \cup S' \cup R \cup \{0, N+1\}$ |
| $V$ | Subset of vertices, i.e. $V = G' \cup D' \cup S'$ |
| $V_0, V_{N+1}, V_{0,N+1}$ | $V_0 = V \cup \{0\}, V_{N+1} = V \cup \{N+1\}, V_{0,N+1} = V \cup \{0, N+1\}$ |
| $\mathcal{A}_c$ | Set of walking arcs from customers' origins to meeting points, *i.e.* $\mathcal{A}_C = \{(r,j) \mid r \in R, j \in G'\}$ |
| $\mathcal{A}_B$ | Set of bus arcs |
| *Parameters and auxiliary variables* | |
| $0, N+1$ | Two duplicate instances of the depot |
| $T$ | Planning horizon |
| $A_i^k$ | Arrival time of bus $k$ at vertex $i$ |
| $B_i^k$ | Beginning time of service of bus $k$ at vertex $i$ |
| $p_i^k$ | Indicator: 1 if node $i$ is visited by bus k, 0 otherwise |
| $W_i^k$ | Waiting time of bus $k$ at node $i$ |
| $q_i^k$ | Passenger load of bus $k$ when leaving vertex $i$ |
| $E_i^k$ | The battery energy level of bus $k$ when arriving at the vertex $i$ |
| $Q^k$ | Capacity of bus $k$ |
| $E_{min}^k, E_{max}^k, E_{init}^k$ | The minimum, maximum, initial state of charge (SOC) of bus $k$ |
| $w_{ri}$ | Walking distance from customer $r$ origin to meeting point $i \in G'$ |
| $c_{ij}$ | Distance from vertex $i$ to vertex $j$ |
| $t_{ij}$ | Bus travel time from vertex $i$ to vertex $j$ $\forall i, j \in V_{0,N+1}$. Note that $t_{rj}$ is the walking time from customer $r$ origin to meeting point $j$, $\forall r \in R, j \in G'$. |
| $L_i$ | Maximum ride time for customers picked up at node $i$. Calculated as 'straight line' ride time multiplied by a detour factor. |
| $w_{max}$ | Maximum walking distance |
| $u_i$ | Service time at vertex $i \in V$ |
| $e_i, l_i$ | Earliest and latest starting times of service at vertex $i \in V_{0,N+1}$ |
| $d_r$ | Drop-off transit station dummy node of customer $r$ |
| $\alpha_s$ | The charging rate of charger $s \in S'$ |
| $\beta^k$ | Energy consumption rate per kilometer traveled for bus $k$ |
| $M$ | Large positive number |
| $\lambda$ | Weighting coefficient for the objective function |
| *Decision variables* | |
| $y_{ri}^k$ | Indicator: 1 if customer $r$ is assigned to bus $k$ and meeting point $i$, 0 otherwise |
| $x_{ij}^k$ | Indicator: 1 if arc $(i,j)$ is traversed by bus $k$, 0 otherwise |
| $\tau_s^k$ | Charging duration for bus $k$ at charger $s$, $s \in S'$ |

**Appendix B. Customer-to-meeting-point-assignment problem.**

For clarity, we provide the non-linear mixed-integer programming formulation of the customer-to-meeting-point problem. Note that it can be easily reformulated as a MILP without difficulty. Note that this problem is the first stage of our metaheuristic algorithm. Its solution serves as input for the subsequent routing routing optimization. The MP-EFCS problem (Eqs. (1)-(34)) is represented on a directed graph with a set of nodes, i.e. $G'$(dummy meeting points), $D'$(dummy transit stations), $S'$(dummy chargers) and $\{0, N+1\}$(two duplicated depots), and two types of arcs, i.e. $\mathcal{A}_C$ and $\mathcal{A}_B$. To solve the problem efficiently, the directed graph is structured in multiple layers, ordered by the departure time of trains (transit service). Each layer $\ell$ is characterized by a pair of transit station ID and the departure time of trains. For each layer $\ell$, there are a subset of dummy meeting points nodes $G'_\ell$, dummy transit station nodes $D'_\ell$, and requests $R_\ell$ aiming to arrive at transit stations for the departure time associated with the layer $\ell$. There are arcs connected between nodes of the same layers and between different layers. Infeasible arcs and inactivated meeting points are trimmed by a pre-processing procedure.

The assignment problem aims to assign customers to nearby dummy meeting points of the same layer to minimize total system cost, measured as a weighted sum of customers' walking time and vehicle routing time between the assigned (activated) meeting points (Eq. (B1)). Note that $\rho$ is an algorithmic parameter to be tuned.

$$\text{Min } \lambda_2 \sum_{r\in R} \sum_{j\in G'_{\ell(r)}} t_{rj}\, y_{rj} + \rho\lambda_1 \sum_{\ell\in\mathcal{L}} \sum_{i\in G'_\ell} \sum_{j\in G'_\ell} t_{ij}\, \theta_i\theta_j \tag{B1}$$

$$c_{rj}y_{rj} \le w_{max}, \forall r \in R, j \in G'_{\ell(r)} \tag{B2}$$

$$\sum_{j\in G'_{\ell(r)}} y_{rj} = 1, \forall r \in R \tag{B3}$$

$$\sum_{r\in R} y_{rj} \le Q_{max}, \forall j \in G' \tag{B4}$$

$$\sum_{r\in R} y_{rj} \le M\theta_j, \forall j \in G' \tag{B5}$$

$$\theta_j \in \{0,1\}, \forall j \in G' \tag{B6}$$

$$y_{rj} \in \{0,1\}, \forall r \in R, j \in G'_{\ell(r)} \tag{B7}$$

Constraint (B2) states that customer (request) $r$ can be assigned to nearby dummy meeting points of the same layer within a maximum walking distance. Constraint (B3) states that each customer needs to be assigned to one meeting point. Constraint (B4) states that the number of assigned customers cannot exceed the vehicle capacity. Constraint (B5) ensures the consistency between the decision variables $y$ and $\theta$. This MILP problem can be solved with little computational time using state-of-the-art commercial MIP solvers. Note that customers with no feasible meeting points within the maximum walking distance are rejected.